\documentclass{article}
\newcommand{\proof}{\noindent {\bf Proof: }}

\newcommand{\corollary}{\noindent {\bf Corollary:}}
\newcommand{\ws}{\hspace{4pt}}
\newtheorem{theorem}{Theorem}

\newtheorem{proposition}{Proposition}
\newtheorem{lemma}{Lemma}
\newtheorem{defi}{Definition}
\usepackage[]{amssymb}
\usepackage[]{amsmath}
\usepackage[]{eucal}

\begin{document}

\title{Weighted Fej\'er Constants and Fekete Sets}
\author{\'A. P. Horv\'ath \footnote{supported by Hungarian National Foundation for Scientific Research, Grant No.  K-100461 \newline{\it Key words}: interpolation, Hermite-Fej\'er, stable and most economical, Fekete sets, Gr\"unwald operator  \newline {\it 2000 MS Classification}: 41A05, 41A36}}\date{}
\maketitle

\begin{abstract}We give the connections among the Fekete sets, the zeros of orthogonal polynomials, $1(w)$-normal point systems, and the nodes of a stable and most economical interpolatory process via the Fej\'er contants. Finally the convergence of a weighted Gr\"unwald interpolation is proved.
\end{abstract}

\section{Introduction}

L. Fej\'er introduced the so-called Hermite-Fej\'er interpolatory process, and in 1934 he gave the definition of normal- and $\varrho$-normal system of nodes for which the Hermite-Fej\'er interpolation is a positive interpolatory process. The surprising nice convergence properties of Lagrange, Hermite and Hermite-Fej\'er operators on $\varrho$-normal systems were proved by L. Fej\'er, G. Gr\"unwald, etc. 
On the other hand the experiences in electrostatics ensure a system of nodes: the Fekete set, which has uniform distribution in some sense, so it must be a good set for interpolation. The system of zeros of orthogonal polynomials has very similar properties, as it it well-known. From another point of view, Egerv\'ary and Tur\'an asked, that is it possible to find an interpolatory process, and a system of nodes together, such that the interpolatory polynomial has the minimal degree, and the operator has the minimal norm. The above-mentioned point systems can be a suitable system of nodes for an interpolatory process in general sense and also with respect to the Egerv\'ary-Tur\'an problem.

The primary aim of this note to revisit the connections among that sets of nodes, and interpolatory problems investigated e.g. in \cite{f}, \cite{h1}, \cite{h2}, \cite{i}, \cite{j}, \cite{rs}. In the next section, we summarize and reformulate these results, and complete them, when the original statement proved only in classical cases. It will be pointed out, that in these equivalences the so-called Fej\'er constants (see(3)) play the key role, that is the characterization of this special system of nodes is ensured by the Fej\'er constants. 

\medskip

As an application of the results of the second section, in the third section we prove a convergence theorem on Gr\"unwald interpolatory process on the real line for Freud-type weights. As it turned out, giving the weighted Fekete sets with respect to a fixed weight is difficult. (However, there are several methods of giving approximating Fekete sets.) The zeros of orthogonal polynomials are Fekete sets for some varying weights. Unfortunately these varying weights tend to zero locally uniformly, so interpolation on Fekete sets in this sense gives only trivial (convergent) processes. The investigation of these weights at infinity leads to define a weighted Gr\"unwald operator (see (11)), which has rather nice convergence properties. Comparing this result with the previous ones of \cite{lu}, \cite{ssz}, it turns out that the convergence is valid here for a wider function class.

\section{Connections}

At first we give the definition of classes of weights in question.
\begin{defi} Let $\Sigma \subset \mathbb{C}$ is a closed set. $w$ is quasi-admissible on $\Sigma$, if it is nonnegative, upper semi-continuous, and if $\Sigma$ is unbounded, $\lim_{|z| \to \infty \atop z \in \Sigma}|z|w(z)=0$. It is admissible, if $\mathrm{cap} \{z\in \Sigma: w(z)>0\}>0$.
Let us call an admissible weight as approximating on $(a,b)\subset \mathbb{R}$, if it has finite moments, it is twice differentiable and $\left(\log\left(\frac{1}{w}\right)\right){''}\geq 0$ on $(a,b)$, and if $a$ is finite, then $\lim_{ x\to a+ }\frac{w(x)}{x-a}=0$, and if $b$ is finite, then $\lim_{ x\to b-} \frac{w(x)}{b-x}=0$.
\end{defi}

\begin{defi}\cite{st}III.1 Let $w$ be a quasi-admissible weight on a closed set $\Sigma \subset \mathbb{C}$. Then $\mathcal{F}_n$ are called n-th weighted Fekete sets associated with $w$, if the supremum below is attained at the set $\mathcal{F}_{n,w}=\{x_1,\dots ,x_n\}$.

$$d_{n,w}=\sup_{z_1, \dots ,z_n \in \Sigma }d_{n,w}(z_1, \dots ,z_n)$$ \begin{equation}=\sup_{z_1, \dots ,z_n \in \Sigma }\left(\prod_{1\leq i<j\leq n}|z_i-z_j|w(z_i)w_(z_j)\right)^{\frac{2}{n(n-1)}}\end{equation}\end{defi}

Usually these points are not unique, but in one dimension by some restrictions on the weight, uniqeness can be proved. In the classical, unweighted case
on $[-1,1]$, the result is proved by Popoviciu (cf. \cite{sze} Ch 6.7 p. 139., and the reference therein). In weighted case, after some restrictions on the weight a representation of Fekete points was given by M. E. H. Ismail (\cite{i} Thms. 2.1, 2.4), wich ensures the unicity of the Fekete sets as well. In the followings the one-dimensional case will be investigated.

\medskip

Now let us deal with the weighted Lagrange interpolatory polynomials on a system of nodes $X=\{x_{k,n}, k=1,\dots ,n; n\in \mathbb{N}\}$. Let $l_k(x)=\frac{\omega(x)}{\omega^{'}(x_k)(x-x_k)},$ where  $\omega(x)=\prod_{k=1}^n(x-x_k)$ (denoting $x_k=x_{k,n}, k=1, \dots ,n$) the fundamental polynomials of the Lagrange interpolation, and let $w(x)=e^{-Q(x)}$ be an  approximating weight. The properties of $L_{k,w,X}(x)=L_{k,w}(x)=w(x)\frac{l_k^2(x)}{w(x_k)}$ will be investigated. It is clear, that $L_{k,w}(x_k)=1$, that is the $\sup$-norm of this weighted polynomial is at least 1. If this $\sup$-norm is equal to one, then $L_{k,w}$ has a maximum at the point $x_k$, that is
$$(L_{k,w})^{'}(x_k)=w(x)\frac{l_k^2(x)}{w(x_k)}\left(\frac{w^{'}(x_k)}{w(x_k)}+\frac{2l_k^{'}(x_k)}{l_k(x_k)}\right)$$ \begin{equation}=w(x)\frac{l_k^2(x)}{w(x_k)}\left(-Q^{'}(x_k)+\frac{\omega^{''}}{\omega^{'}}(x_k)\right)=0,\end{equation}
which ensures that
\begin{equation}C_{k,w}:=C_{k,w,X}=\frac{\omega^{''}}{\omega^{'}}(x_k)+\frac{w^{'}}{w}(x_k)=0.\end{equation}
This is the case, when $X$ is a Fekete set with respect to $w^{\frac{1}{2(n-1)}}$, namely
$$L_{k,w,X}=\frac{\prod_{1\leq l \leq n \atop l\ne k}\left((x-x_l)^2w^{\frac{2}{2(n-1)}}(x)w^{\frac{2}{2(n-1)}}(x_l)\right)}{\prod_{1\leq l \leq n \atop l\ne k}\left((x_k-x_l)^2w^{\frac{2}{2(n-1)}}(x_k)w^{\frac{2}{2(n-1)}}(x_l)\right)}$$ $$\times\frac{\prod_{1\leq i<j \leq n \atop i,j \ne l}\left((x_i-x_j)^2w^{\frac{2}{2(n-1)}}(x_i)w^{\frac{2}{2(n-1)}}(x_j)\right)}{\prod_{1\leq i<j \leq n \atop i,j \ne l}\left((x_i-x_j)^2w^{\frac{2}{2(n-1)}}(x_i)w^{\frac{2}{2(n-1)}}(x_j)\right)}\leq 1,$$
because in the denominator appears $d_{n,w^{\frac{1}{2(n-1)}}}^{n(n-1)}$.

\medskip

It will turn out in the followings, that the behavior of the constants $C_{k,w}$ as an indicator, shows the properties of the point systems, interpolatory systems and operators. Emphasizing the importence of these constants, let us call them as "{\bf Fej\'er constants}".

\medskip

Following carefully the proof of the above mentioned theorem of Ismail (\cite{i},Thm. 2.1), we get the following

\begin{proposition}Let $w$ be an approximating weight on an interval $(a,b)$. Then $d_{n,w^{\frac{1}{2(n-1)}}}^{n(n-1)}(z_1, \dots ,z_n)$ attains its maximum on $(a,b)$ at a unique set $\mathcal{F}_{n,w^{\frac{1}{2(n-1}}}$, for which the following characterization is valid.
\begin{equation} \mathcal{F}_{n,w^{\frac{1}{2(n-1)}}}=\{x_1, \dots ,x_n\} \ws \ws \mbox{if and only if} \ws \ws C_{k,w}=0,\ws \ws k=1, \dots , n. \end{equation}\end{proposition}

At first we have to note here, that finite moments are not necessary in this statement. According to Ismail \cite{i}, the proof of this theorem is the following: taking the partial derivatives of $\log d_{n,w^{\frac{1}{2(n-1)}}}^{n(n-1)}$, it turns out, that \\
$\frac{\partial}{\partial x_j}\log d_{n,w^{\frac{1}{2(n-1)}}}^{n(n-1)}(x_1, \dots ,x_n) =0$ $j=1, \dots , n$, if and only if $C_{k,w}=0,\ws \ws k=1, \dots , n$. Computing the Hessian, it can be seen, that $-H$ is always positive definite, so recalling the boundary condition on $w$, we get that the maximum-set is unique, that is it is the unique solution of the equation system: $C_{k,w}=0,\ws \ws k=1, \dots , n$. Independently of the previous chain of ideas, an elementary proof on unicity can be given.

\begin{proposition}Let $w$ be an admissible, continuous weight on $\mathbb{R}$ such that $\log\frac{1}{w}$ is convex. Then the associated weighted Fekete sets are unique.\end{proposition}

\proof
Contrary, let $\{x_i\}_{i=1}^n$ and $\{y_i\}_{i=1}^n$ are Fekete points with respect to $w$ enumerated in increasing order, and let $z_i=\frac{x_i+y_i}{2}$. Then because of the ordering of the points, and the log-convexity of the weight, by the arithmetic-geometric mean inequality
$$|z_i-z_j|w(z_i)w(z_j)=\left|\frac{(x_i-x_j)+(y_i-y_j)}{2}\right|w\left(\frac{x_i+y_i}{2}\right)w\left(\frac{x_j+y_j}{2}\right)$$ $$=\frac{|x_i-x_j|+|y_i-y_j|}{2}w\left(\frac{x_i+y_i}{2}\right)w\left(\frac{x_j+y_j}{2}\right)$$ $$\geq \sqrt{|x_i-x_j|}\sqrt{|y_i-y_j|}\sqrt{w(x_i)}\sqrt{w(y_i)}\sqrt{w(x_j)}\sqrt{w(y_j)},$$
where the inequality is an equality if and only if $x_i=y_i$ for all indices, wich establishes the uniqueness.

\medskip

For special weights, the Fekete sets are the zeros of some orthogonal polynomials (cf. \cite{i}, \cite{h2}). Before setting the precise statement we need some definitions.

\begin{defi} Let $w=e^{-Q}$ be an approximating weight on $(a,b)$. Let
\begin{equation}A_n(x)=\varrho_n\int_a^bp_{n,w}^2(t)w(t)\frac{Q^{'}(t)-Q^{'}(x)}{t-x}dt, \end{equation}
where $p_{n,w}=\gamma_nx^n+\dots$ is the $n^{th}$ orthonormal polynomial with respect to $w$, and $\varrho_n=\frac{\gamma_{n-1}}{\gamma_n}$ \end{defi}
Now we can define our weights:
\begin{defi}Let $w$ be as in the previous definition.
\begin{equation}w_n(x)=\frac{w(x)\varrho_n}{A_n(x)}\end{equation}\end{defi}

In the following investigations the constant $\varrho_n$ has not any role, but it will come into the picture inconnection with a convergence theorem in the next section. Let us see some examples on $\frac{A_n(x)}{\varrho_n}$ (\cite{i}), which in classical cases are different only in normalization from the weights $w_1$ (\cite{rs}), for which the derivatives of $p_{n,w}$-s are orthogonal :

\noindent{\bf Example:} 

\noindent (1) If $w=e^{-x^2}$, 
$$\frac{A_n(x)}{\varrho_n}=2,$$
that is $w_n=\frac{1}{2}w$ independently of $n$ and $x$, and here $w_1=w=2w_n$.

\medskip

\noindent (2) If $w=x^{\alpha}e^{-x}$,
$$\frac{A_n(x)}{\varrho_n}=\frac{1}{x},$$
that is $w_n=x^{\alpha+1}e^{-x}=xw$ independently of $n$, and here $w_1=w_n$.

\medskip

\noindent (3) If $w=(1-x)^{\alpha}(1+x)^{\beta}$,
$$\frac{A_n(x)}{\varrho_n}=\frac{\alpha+\beta+1+2n}{1-x^2},$$
that is $w_n=\frac{1}{\alpha+\beta+1+2n}(1-x)^{\alpha+1}(1+x)^{\beta+1}$, and here $w_1=(\alpha+\beta+1+2n)w_n$.

\medskip

\noindent (4) If $w=e^{-x^4}$, 
$$\frac{A_n(x)}{\varrho_n}=2(x^2+\varrho_n^2+\varrho_{n+1}^2),$$
that is $w_n=\frac{1}{2(x^2+\varrho_n^2+\varrho_{n+1}^2)}w$.

\medskip

From another point of view $w_n$ has also an importance. Denoting by $p_n\sqrt{w_n}=z_n$, it satisfies the following differential equation with some $\Phi_n$ (cf. \cite{mh}, Th. 3.6.):
\begin{equation}z_n^{''}(x)+\Phi_n(x)z_n(x)=0\end{equation}

\medskip 

In the next statement we reformulate the results of Ismail, Rutka and Smarzewski (cf. \cite{i}, \cite{rs}).

\begin{proposition} Let $w_n$ be as in the definitions above, and let us assume that $w_n$ is an approximating weight. Then
\begin{equation} C_{k,w_n}=0,\ws \ws k=1, \dots , n \ws\ws\mbox{if and only if} \ws \ws \{x_k\} \ws\ws \mbox{the zeros of} \ws\ws p_{n,w}\end{equation}\end{proposition}

The proof of this statement depends on the differential equation of orthogonal polynomials. The equation system on $C_k$-s means that the differential equation fulfils at the points $x_k, \ws k=1, \dots , n $. In the classical cases, it is a Sturm-Liouville equation, that is there are polynomials of degree $n$ in the differential equation, which is realized at $n$ points. In general cases unicity is used.

\medskip

Normal and $\varrho$-normal point systems were introduced on $[-1,1]$ by L. Fej\'er in 1934 (\cite{f}). The weighted analogon of this definition was given in \cite{h1}. The original aim of these definitions was assuring the positivity of the Hermite-Fej\'er interpolatory operator. The limit case, when $\varrho=1$ was investigated on the weighted real line in \cite{h2}. Here this last definition is cited only.

\begin{defi} Let $w$ be an approximating weight on $(a,b)$. A system of nodes $X=\{x_{k,n}, k=1,\dots ,n; n\in \mathbb{N}\}$ is $1(w)$-normal, if there is an $L>1$ such that
\begin{equation}|x_{k,n}|<La_n,\end{equation}
where $a_n$ is the M-R-S number, and
\begin{equation}w(x)\sum_{k=1}^n\frac{l_k^2(x)}{w(x_k)} \leq 1, \ws\ws x\in \mathbb{R},\end{equation}
where $l_k(x)$-s are the fundamental polynomials of the Lagrange interpolation.\end{defi}

In this definition the kernel function of the Gr\"unwald operator appears. Here we will follow the notations of \cite{lu} and \cite{ssz}, that is the weighted Gr\"unwald operator on the nodes $\{x_k\}_{k=1}^n$ with respect to an $f$ is
\begin{equation}w(x)Y_n(f,x)=w(x)\sum_{k=1}^n l_k^2(x)f(x_k)\end{equation}
Mostly the boundedness of the operator-norm ensures the convergence of the interpolatory process. The boundedness by one, is a very special criterium. This is the case for instance, when the reciprocal of the weight function has non-negative even derivatives, and the Gr\"unwald operator coincides with the Hermite-Fej\'er one. Also on this chain of ideas the Fej\'er constants play the key role. More precisely, with the notations above, the weighted Hermite interpolatory polinomial (with some weight $w$) of a differentiable function can be expressed as (cf. \cite{h2})
$$w(x)H_n(f,f^{'},x)=w(x)\sum_{k=1}^n\frac{(1-C_{k,w}(x-x_k))l_k^2(x)}{w(x_k)}(fw)(x_k)$$ \begin{equation}+w(x)\sum_{k=1}^n\frac{(x-x_k)l_k^2(x)}{w(x_k)}(fw)^{'}(x_k),\end{equation}
and the corresponding weighted Hermite-Fej\'er operator is
\begin{equation}w(x)H_{n,w}(f,x)=w(x)\sum_{k=1}^n\frac{(1-C_{k,w}(x-x_k))l_k^2(x)}{w(x_k)}(fw)(x_k),\end{equation}
which coincides with the weighted function at the nodes $\{x_k\}_{k=1}^n$, and which has zero derivatives at the nodes. Furthermore by the definition of the Fej\'er constants, $H_{n,w}(f,x)$ is the (unweighted) Hermite interpolatory polynomial of $\frac{1}{w}$. So when the Fej\'er constants are zero
$$Y_{n,w}(x):=w(x)Y_n\left(\frac{1}{w},x\right)=w(x)\sum_{k=1}^n\frac{l_k^2(x)}{w(x_k)}=w(x)H_{n,w}(\frac{1}{w},x)$$ \begin{equation}=w(x)\sum_{k=1}^n\frac{(1-C_{k,w}(x-x_k))l_k^2(x)}{w(x_k)}=w(x)H_n\left(\frac{1}{w},\left(\frac{1}{w}\right)^{'},x\right)\end{equation}
is the Hermite interpolatory polynomial of $\frac{1}{w}$ with respect to the nodes: $\{x_k\}_{k=1}^n$. So the following connections are established.

\begin{proposition} Let $w$ be a weight as above. 

If a system of nodes $\{x_{k}\}$ is $1(w)$-normal, then $C_{k,w}=0,\ws \ws k=1, \dots , n $. 

On the other hand, let us suppose further, that $\left(\frac{1}{w}\right)^{2n}\geq 0$ on $|x|\leq La_n$. Now 

if $C_{k,w}=0,\ws \ws k=1, \dots , n $ then the system of nodes is $1(w)$-normal.\end{proposition}

\proof
If $\{x_{k}\}$ is $1(w)$-normal, then $w(x)\frac{l_k^2(x)}{w(x_k)}\leq 1, \ws k=1, \dots , n$ (see (10)), so $C_{k,w}=0,\ws \ws k=1, \dots , n $. According to (14), by
the error formula of the Hermite interpolation, it is clear, that $1-w(x)\sum_{k=1}^n\frac{l_k^2(x)}{w(x_k)}\geq 0$, when $\left(\frac{1}{w}\right)^{(2n)} \geq 0$ on $|x|\leq La_n$.

\medskip

The Egerv\'ary-Tur\'an interpolatory problem (cf. \cite{rs}, and the references therein) is to find an interpolatory process of lowest degree, and of smallest norm.
Below we denote by $\hat{l}_k(x)$ each polynomial of arbitrary degree for which $\hat{l}_k(x_i)=\delta_{ki}, \ws i=1, \dots , n$.
\begin{defi} Let $w$ be as in Definition 5. The interpolatory system of polynomials $\hat{l}_k(x), \ws k=1, \dots , n$ is $w$-stable on $(a,b)$ if for all $y_1, \dots , y_n \geq 0$
\begin{equation} 0\leq w(x)\sum_{k=1}^n\frac{\hat{l}_k(x)}{w(x_k)}y_k \leq \max_ky_k, \ws\ws x\in (a,b).\end{equation}
A $w$-stable interpolatory system on $(a,b)$ is most economical, if 
\begin{equation} \sum_{k=1}^n \deg \left(\hat{l}_k(x)\right)\end{equation}
is minimal.\end{defi}

Let us remark that if the weight function tends to zero quickly at the boundary points of the fundamental interval, then the $w$-stability of the Gr\"unwald operator  coincides with the $1(w)$-normality of the nodes. It is proved for all the classical weights (cf. \cite{rs}, Thm. 2.3), that an interpolatory system is $w_n$-stable and most economical, if and only if it is the Gr\"unwald operator on the zeros of $p_{n,w}$. From the previous investigations, similarly to the classical cases, we can state the parallel theorem for general weights. 

Let us denote by

$$I_{n,w}(x):=w(x)\sum_{k=1}^n\frac{\hat{l}_k(x)}{w(x_k)}$$

\begin{proposition}  Let $w$ be an approximating weight on an interval $(a,b)$. 

\noindent If $I_n(x)$ is $w$-stable and most economical, then
\begin{equation}C_{k,w}=0,\ws \ws k=1, \dots , n.\end{equation} 

\noindent Let us assume further that $\left(\frac{1}{w}\right)^{(2n)} \geq 0$ on $|x|\leq La_n$. 

\noindent If $C_{k,w}=0,\ws \ws k=1, \dots , n ,$ then

\begin{equation} I_n(x)=Y_{n,w}(x)\ws\ws\mbox{is $w$-stable and most economical}\end{equation}

\end{proposition}

\proof
As it was pointed out eg. in \cite{rs},  if an interpolatory process $I_n(x)$ is $w$-stable and most economical, it must be the Gr\"unwald operator, because by the positivity of the operator, $\hat{l}_k(x)$ has zeros at the points $x_i, i=1,\dots ,n, i\neq k$ of even multiplicity, that is $\sum_{k=1}^n \deg \left(\hat{l}_k(x)\right)\geq 2n(n-1)$. It is realized by $Y_{n,w}$. As it was shown in Statement 4, if $Y_{n,w}$ has maxima at $x_k$-s then $C_k$-s are zero. The opposite direction is also follows from Statement 4.

\medskip

Finally enumerating the properties discussed above, we can summarize these results as it follows.

$$\begin{array}{ll} \mbox{(\bf{A})} \ws\ws\ws C_{k,w}=0,\ws \ws k=1, \dots , n \ws\ws\ws (\mathrm{\bf{A}^{'}}) \ws\ws\ws C_{k,w_n}=0,\ws \ws k=1, \dots , n  \\
\mbox{(\bf{B})}\ws\ws\ws \mathcal{F}_{n,w^{\frac{1}{2(n-1)}}}=\{x_1, \dots ,x_n\}  \ws\ws\ws (\mathrm{\bf{B}^{'}}) \ws\ws\ws \mathcal{F}_{n,w_n^{\frac{1}{2(n-1)}}}=\{x_1, \dots ,x_n\}\\
 \mbox{(\bf{C})}\ws\ws\ws p_{n,w}(x_k)=0,\ws\ws  k=1, \dots , n\\
 \mbox{(\bf{D})} \ws\ws\ws  Y_{n,w}(x) \ws\ws\mbox{is $w$-stable and most economical}\\ (\mathrm{\bf{D}^{'}}) \ws\ws\ws Y_{n,w_n}(x) \ws\ws\mbox{is $w_n$-stable and most economical}\\
  \mbox{(\bf{E})} \ws\ws\ws \{x_1, \dots ,x_n\} \ws\ws\mbox{is $1(w)$-normal}\ws\ws\ws (\mathrm{\bf{E}^{'}}) \ws\ws\ws \{x_1, \dots ,x_n\} \ws\ws\mbox{is $1(w_n)$-normal}\end{array} $$
Through the equivalence of all the above mentioned properties with property ({\bf A}) (or $(\mathrm{\bf{A}^{'}})$), that is the Fej\'er constants are zero, one can get

\corollary

Let $w$ be an admissible, approximating weight on an interval $(a,b)$. If $\left(\frac{1}{w}\right)^{(2n)} \geq 0$ on $(a,b)$, then ({\bf A}),({\bf B}),({\bf D}),({\bf E}) are equivalent, and if
$\left(\frac{1}{w_n}\right)^{(2n)} \geq 0$ on $(a,b)$, then $(\mathrm{\bf{A}^{'}}), (\mathrm{\bf{B}^{'}}),$ ({\bf C}), $(\mathrm{\bf{D}^{'}}), (\mathrm{\bf{E}^{'}})$ are equivalent.
 
 \medskip 
 
We have to show an example on the second assumption .

\noindent{\bf Example:}
Let 
\begin{equation} Q(x)=\sum_{k=0}^m d_k x^{2k}, \hspace{1cm} d_k \geq 0, \ws k=1, \dots ,m, \end{equation}
and let $w(x)=e^{-Q(x)}$. For these special Freud-type weights $\left(\frac{1}{w_n}\right)^{(2n)} \geq 0$ on $\mathbb{R}$ for all $n \in \mathbb{N}$ . 
According to the Leibniz rule it is enough to show that $\left(\frac{A_n}{\varrho_n}\right)^{(j)}\left(\frac{1}{w}\right)^{(2n-j)}>0$ for $j=1, \dots ,2n$. Because 
$$\frac{\partial^j\left(\bar{Q}(t,x)\right)}{\partial x^j}= \sum_{k=1}^m 2kd_k\frac{\partial^j\left(\frac{t^{2k-1}- x^{2k-1}}{t-x}\right)}{\partial x^j}=\sum_{k=\lceil\frac{j}{2}\rceil+1}^m2kd_k\sum_{l=j}^{2k-2}b_lt^{2k-2-l}x^{l-j},$$
where $b_l$-s are positive, taking into consideration that $w$ is an even weight function, (and so $p_n^2(w)$ is also even), one can see that $$\left(\frac{A_n}{\varrho_n}\right)^{(j)}=\sum_{k=\lceil\frac{j}{2}\rceil+1}^m 2kd_k\sum_{l=j}^{2k-2}b_l\int_{\mathbb{R}}p_n^2(w,t)w(t)t^{2k-2-l}dtx^{l-j}$$ 
is a polynomial of $x$ with nonnegative coefficients, and all the exponents of this polynomial are even if $j$ is even and are odd if $j$ is odd. By a simple induction one can see that 
$$\left(\frac{1}{w(x)}\right)^{(j)} = p(j,x)e^{Q(x)},$$
where $p(j,x)$ is a polynomial having the same properties as the previous one. Because $j$ and $2n-j$ have the same parity, 
$$\left(\frac{1}{w_n}\right)^{(2n)}(x)=p(x)e^{Q(x)},$$
where $p(x)$ is a polynomial with even exponents and positive coefficients, so it is positive on the real line for all $n\in \mathbb{N}$.

\medskip

Finally we have to remark that the assumption $\left(\frac{1}{w}\right)^{(2n)} \geq 0$ seems to be assymetric, and it is necessary only because of the method of the proof by Hermite interpolation. The question that can it be weakened or not, is unsolved yet.

\section{Interpolation}

In this section, let $w=e^{-Q}$ be a three times continuously differentiable Freud weight on $\mathbb{R}$, that is we suppose that $Q$ is even, $Q^{'}>0$ on $(0,\infty)$, and for some $A,B\geq 2$;  $A\leq\frac{(xQ^{'}(x))^{'}}{Q^{.}(x)}\leq B$ on $(0,\infty)$, moreover there is a constant $c$ such that for every $|x|\geq 1$, $\left|\frac{xQ^{(3)}(x)}{Q^{''}(x)}\right|\leq c$. By these assumptions there is a $d \geq 1$ such that $Q^{''}(x) \geq \frac{1-(B-1)^2-c(B-1)}{x^2}$, when $|x|\geq d$. Now we can define

\begin{defi} With $d>1$ given above, let
\begin{equation} \tilde{w}(x)=\left\{ \begin{array}{ll} w(x), |x| \leq 1\\
w(x)\frac{x}{Q^{'}(x)}, |x| \geq d\\

\mbox{twice continuously differentiable}, \mbox{elsewhere}\end{array}\right.\end{equation}

Furthermore we assume that $\log\frac{1}{\tilde{w}}$ has positive and continuous first and second derivatives on $(0,\infty)$.
\end{defi}

Let us remark at first that $Q^{'}(1)\leq Q^{'}(d)+\frac{d}{Q^{'}(d)}\left(\frac{Q^{'}(x)}{x}\right)^{'}(d)$, because $Q^{'}$ is increasing, and the second member of the right-hand side is positive, when $A\geq 2$. That is a suitable connection can be defined between the two parts of $\log\frac{1}{\tilde{w}}$. 

As  usually we define
\begin{defi}  $$C_{\tilde{w}}=\{f\in C(\mathbb{R})|\lim_{|x| \to \infty} (f\tilde{w})(x) =0$$\end{defi}

Let $Y_n(f,\cdot)$ be as in (10), the Gr\"unwald operator on the zeros of $p_{n,w}$. Now we have the following

\begin{theorem} Let $f \in C_{\tilde{w}}$ Then
\begin{equation}\lim_{n\to \infty}\|(Y_n(f)-f)\tilde{w}\|=0\end{equation}\end{theorem}

Comparing this theorem with Cor. 2. of \cite{ssz}, we can see, that we have two different weights in this theorem, but when $A\geq 2$, then the function class is wider here, that is the fuctions can grow more quickly at infinity.

The previous definition of the weight was inspirated by the next lemma. Investigating the weights $w_n$ from the previous section, it turns out, that however $w_n$ tends to zero locally uniformly when $n$ tends to infinity, the behavior of $w_n$-s are the same at infinity. It means, that the Gr\"unwald operator on Fekete points with respect to the varying weights $w_n$ has trivial convergence properties, but it allows to find a non-trivial process, as it is given in the theorem.

The following estimation of $A_n$ is valid.

\begin{lemma} Let $w$ be as above, and let $A\geq 2$. Let $L_0$ be a constant such that $\frac{L_0}{2}a_n>a_{2n+[A-1]+1}$. For every $L>L_0$
\begin{equation}\frac{A_n(x)}{\varrho_n}\sim\left\{\begin{array}{ll}\frac{n}{a_n^2}, \ws\mbox{if}\ws |x| \leq La_n\\
\frac{Q^{'}(x)}{x},\ws\mbox{if}\ws |x| \geq La_n\end{array}\right.,\end{equation}
where the constants in $"\sim"$ depend only on $L$, but they are independent of $n$.\end{lemma}

\proof

At first we have to note that such an $L_0$ exists by \cite{lelu} 5.9.
The first line of the inequality is proved by H. N. Mhaskar (\cite{mh}, Prop. 3.7.). To prove the second line we have to divide the integral to some parts. Since $A_n$ is even we can choose $x>La_n$.
$$\frac{A_n(x)}{\varrho_n}=\int_{|t|\leq \frac{x}{2}}p_n^2(t)w(t)\frac{Q^{'}(t)-Q^{'}(x)}{t-x}dt+\int_{|t|>\frac{x}{2}}(\cdot)dt=I_1+I_2$$
If $|t|\leq  \frac{x}{2}$, using that $Q$ is convex on $\mathbb{R}$, and estimating the denominator as $\frac{x}{2} \leq |t-x| \leq \frac{3}{2}x$,
$$\frac{2}{3}\frac{\left(2^{A-1}-1\right)}{2^{A-1}}\frac{Q^{'}(x)}{x}\leq \frac{2}{3}\frac{Q^{'}(x)-Q^{'}\left(\frac{x}{2}\right)}{x}\leq \frac{Q^{'}(t)-Q^{'}(x)}{t-x}$$ $$\leq 2\frac{|Q^{'}(t)|+Q^{'}(x)}{x}\leq 4\frac{Q^{'}(x)}{x},$$
where in the first inequality we used the properties of $Q^{'}$, cf \cite{lelu} 5.3.
 So 
\begin{equation}I_1 \sim\frac{Q^{'}(x)}{x}\int_{|t|\leq \frac{x}{2}}p_n^2w\sim\frac{Q^{'}(x)}{x}\int_{\mathbb{R}}p_n^2w\sim\frac{Q^{'}(x)}{x},\end{equation}
 where in the second "tilde", the lower estimation fulfils by \cite{mh},(2.6), say.

For  $|t|>\frac{L}{2}a_n>4a_n$
\begin{equation}I_2=\int_{|t|>\frac{x}{2}}p_n^2(t)w(t)\frac{Q^{'}(t)-Q^{'}(x)}{t-x}dt
=\int_{|t|>\frac{x}{2}\atop |x-t|\leq 1}(\cdot)+\int_{|t|>\frac{x}{2}\atop |x-t|>1}(\cdot) =I_3+I_4\end{equation}
By the properties of Freud weights, and by \cite{mh},(2.6),
$$I_3 =\int_{|t|>\frac{x}{2}\atop|x-t|\leq 1}p_n^2(t)w(t)Q^{''}(\xi(x,t))dt$$ 
\begin{equation}\leq c \frac{Q^{'}(x)}{x}\int_{ |t|>\frac{x}{2}}p_n^2(t)w(t)dt\leq c_1e^{-c_2n}\frac{Q^{'}(x)}{x}\end{equation}
$$I_4=\int_{-\infty}^{-\frac{x}{2}}p_n^2(t)w(t)\frac{Q^{'}(|t|)+Q^{'}(x)}{|t|+x}dt+\int_{\frac{x}{2}}^{x-1}p_n^2(t)w(t)\frac{Q^{'}(t)-Q^{'}(x)}{t-x}dt$$ \begin{equation}+\int_{x+1}^{2x}p_n^2(t)w(t)\frac{Q^{'}(t)-Q^{'}(x)}{t-x}dt+\int_{2x}^{\infty}(\cdot)=I_5+I_6+I_7+I_8\end{equation}
$$I_5\leq \int_{-\infty}^{-\frac{x}{2}}p_n^2(t)w(t)\frac{Q^{'}(|t|)}{|t|}dt+\int_{-\infty}^{-\frac{x}{2}}p_n^2(t)w(t)\frac{Q^{'}(x)}{x}dt$$ According to \cite{lelu}, 5.2,
$$I_5\leq \int_{\frac{x}{2}}^{\infty}p_n^2(t)w(t)t^{[A-1]+1}dt+c_1e^{-c_2n}\frac{Q^{'}(x)}{x}$$
According to \cite{mh} 2.7,
$$\int_{\frac{x}{2}}^{\infty}p_n^2(t)w(t)t^{[A-1]+1}dt\leq c_1e^{-c_2n}\int_{|t|\leq a_{2n+[A-1]+1}}p_n^2(t)w(t)t^{[A-1]+1}dt$$
Because $A\geq 2$, $\frac{Q^{'}(t)}{t}$ is increasing, so by \cite{lelu} 5.2 and 5.9
$$\int_{|t|\leq a_{2n+[A-1]+1}}p_n^2(t)w(t)t^{[A-1]+1}dt\leq \frac{Q^{'}(x)}{x}\int_{|t|\leq a_{2n+[A-1]+1}}p_n^2(t)w(t)t^2dt$$ $$\leq c a_n^2\frac{Q^{'}(x)}{x}$$
that is 
\begin{equation}I_5\leq c_3e^{-c_2n}\frac{Q^{'}(x)}{x}\end{equation}
If $t\in \left(\frac{x}{2},x-1\right)$ we can write that $x=\lambda t$, where $1<\lambda\leq 2$, that is according to \cite{lelu}, 5.3, recalling that $B\geq A\geq 2$, $\frac{Q^{'}(x)-Q^{'}(t)}{x-t}\leq \frac{Q^{'}(t)}{t}\frac{\lambda^{B-1}-1}{\lambda-1}\leq c(B) \frac{Q^{'}(t)}{t}\leq c(B) \frac{Q^{'}(x)}{x}$, where in the last step we used that $\frac{Q^{'}(t)}{t}$ is increasing. So
\begin{equation}I_6 \leq c \frac{Q^{'}(x)}{x}\int_{\frac{x}{2}}^{x-1}p_n^2(t)w(t)dt \leq c_1e^{-c_2n}\frac{Q^{'}(x)}{x}\end{equation}
As in the previous case, when $x+1<t<2x$, $\frac{Q^{'}(t)-Q^{'}(x)}{t-x}\leq \frac{Q^{'}(x)}{x}\frac{\lambda^{B-1}-1}{\lambda-1}$ so
\begin{equation}I_7 \leq c(B)\frac{Q^{'}(x)}{x} \int_{x+1}^{2x}p_n^2(t)w(t)dt\leq c_1e^{-c_2n}\frac{Q^{'}(x)}{x}\end{equation}
Because $t-x>\frac{t}{2}$ in $I_8$,
$$I_8\leq 2 \int_{2x}^{\infty}p_n^2(t)w(t)\frac{Q^{'}(t)}{t}dt \leq 2\int_{2x}^{\infty}p_n^2(t)w(t)t^{[A-1]+1}dt$$
So as in $I_5$,
\begin{equation}\leq c_1e^{-c_2n}\int_{|t|\leq a_{2n+[A-1]+1}}p_n^2(t)w(t)t^{[A-1]+1}dt\leq c_1e^{-c_2n}\frac{Q^{'}(x)}{x}\end{equation}
That is $I_1 \sim \frac{Q^{'}(x)}{x}$, and $I_2\leq c_1e^{-c_2n}\frac{Q^{'}(x)}{x},$ which proves the second line of the lemma.

\medskip

For the proof of the Theorem we need the following

\begin{lemma}
If $f \in C_{\tilde{w}}$, then by the notation of (13)
\begin{equation} \|Y_{n,\tilde{w}}\|=O(1)\end{equation}\end{lemma}

\proof

At first let $\frac{a_n}{2}\leq |x|\leq La_n$. Here, by \cite{lelu} 5.5, $\frac{x}{Q^{'}(x)}\sim \frac{a_n^2}{n}$. Using that $\frac{Q^{'}(x)}{x}$ is even, and it is increasing on $\mathbb{R}_+$ when $A\geq 2$, we have 
$$Y_{n,\tilde{w}}(x)=\frac{x}{Q^{'}(x)}w(x)\sum_{k=1}^n\frac{l_k^2(x)}{w(x_k)}\frac{Q^{'}(x_k)}{x_k}$$ \begin{equation}\leq c \frac{a_n^2}{n}w(x)\sum_{k=1}^n\frac{l_k^2(x)}{w(x_k)}\frac{Q^{'}(a_n)}{a_n}\leq c w(x)\sum_{k=1}^n\frac{l_k^2(x)}{w(x_k)}=O(1), \end{equation}
For the last equality cf. \cite{sz} (39).

When $|x|>La_n>|x_k|$, then $A\geq 2$ yields that $\frac{x}{Q^{'}(x)}\frac{Q^{'}(x_k)}{x_k}\leq 1$. That is  
\begin{equation}Y_{n,\tilde{w}}(x)\leq c w(x)\sum_{k=1}^n\frac{l_k^2(x)}{w(x_k)}=O(1), \end{equation}
as in the previous case.

Let $d<|x|<\frac{a_n}{2}$.
$$Y_{n,\tilde{w}}(x)=\tilde{w}(x)\sum_{k \atop |x_k|\leq 2|x|}(\cdot)+\tilde{w}(x)\sum_{k \atop |x_k|> 2|x|}(\cdot)=\Sigma_1+\Sigma_2$$
Because as in previously, when $|x_k|\leq 2|x|$, $\frac{x}{Q^{'}(x)}\frac{Q^{'}(x_k)}{x_k}<c(B)$'
\begin{equation}\Sigma_1 = O(1)\end{equation}
$$\Sigma_2 = w(x)p_n^2(x)\frac{x}{Q^{'}(x)}\sum_{k \atop |x_k|> 2|x|}\frac{Q^{'}(x_k)}{w_(x_k)x_k(x-x_k)^2p_n^{'2}(x_k)}$$
Since $\frac{1}{p_n^{'2}(x_k)w(x_k)}\sim \frac{a_n^2}{n}\Delta x_k$, (cf \cite{cr} 4.11, 4.17)
$$\Sigma_2\leq c w(x)p_n^2(x)\frac{x}{Q^{'}(x)}\frac{a_n^2}{n}\sum_{k \atop |x_k|> 2|x|}\frac{Q^{'}(x_k)}{x_k(x-x_k)^2}\Delta x_k$$ 
$$ \leq cw(x)p_n^2(x)\frac{x}{x^2Q^{'}(x)}\frac{a_n^2}{n}\int_{2d}^{a_n}\frac{Q^{'}(x)}{x}dx \leq c \frac{a_n}{n}\frac{1}{xQ^{'}(x)}$$ \begin{equation}\times \left(\frac{n}{a_n^2}+\int_{2d}^{a_n}\frac{Q(x)}{x^2}dx\right) 
\leq c\frac{a_n}{n}\frac{1}{xQ^{'}(x)}Q(a_n)\leq c\frac{1}{dQ^{'}(d)}=O(1)\end{equation}
Here we used that $w(x)p_n^2(x)\leq \frac{c}{a_n}$ for $|x| \leq \frac{a_n}{2}$, cf. \cite{cr} 4.6.

Finally let $|x|\leq d$. Let us remark, that $\frac{\tilde{w}(x)}{w(x)}$ is between two constants on $[1,d]$.
$$Y_{n,\tilde{w}}(x) \leq c w(x)\sum_{k, x_k<2d}\frac{l_k^2(x)}{w(x_k)}+ cw(x)\sum_{k, x_k\geq 2d}\frac{l_k^2(x)}{w(x_k)}\frac{Q^{'}(x_k)}{x_k}=\Sigma_3+\Sigma_4$$
As previously, 
\begin{equation}\Sigma_3 = O(1)\end{equation}
Similarly to the estimation of $\Sigma_2$,
\begin{equation}\Sigma_4\leq c\frac{a_n}{n}\sum_{k, x_k\geq 2d}\frac{Q^{'}(x_k)}{x_k(x-x_k)^2}\Delta x_k \leq \frac{c}{d^2}\frac{a_n}{n}\left(\frac{n}{a_n^2}+Q(a_n)\right)=O(1)\end{equation}

\medskip 

{\bf Proof (of the Theorem):}
Because polynomials are obviously in $C_{\tilde{w}}$, according to the Banach-Steinhaus theorem, the previous lemma ensures the result.

\medskip
\noindent \small{Department of Analysis, \newline
Budapest University of Technology and Economics}\newline
\small{ g.horvath.agota@renyi.mta.hu}

\end{document}